\documentclass[12pt]{article}
\usepackage{amsmath,amsfonts,amssymb,rotating,colordvi}
\usepackage{colordvi}
\usepackage{mathrsfs}
\def\qed{\nopagebreak\hfill{\rule{4pt}{7pt}}}
\def\proof{\noindent {\it{Proof.} \hskip 2pt}}

\newtheorem{theo}{Theorem}[section]

\newtheorem{prop}[theo]{Proposition}

\makeatletter \@addtoreset{equation}{section}
\makeatletter
\def\Eqlfill@{\arrowfill@\Relbar\Relbar\Relbar}
\newcommand{\extendEql}[1][]{\ext@arrow 0359\Eqlfill@{#1}}
\makeatother

\numberwithin{equation}{section}

\newdimen\Squaresize \Squaresize=11pt
\newdimen\Thickness \Thickness=0.7pt
\def\Square#1{\hbox{\vrule width \Thickness
   \vbox to \Squaresize{\hrule height \Thickness\vss
    \hbox to \Squaresize{\hss#1\hss}
   \vss\hrule height\Thickness}
\unskip\vrule width \Thickness} \kern-\Thickness}

\def\Vsquare#1{\vbox{\Square{$#1$}}\kern-\Thickness}

\def\moins{\raise 1pt\hbox{{$\scriptstyle -$}}}
\begin{document}

\begin{center}
{\large \bf  On the $q$-log-convexity conjecture of Sun}
\end{center}

\begin{center}
Donna Q. J. Dou$^{1}$ and Anne X. Y. Ren$^{2}$\\[9pt]

$^1$School of Mathematics, Jilin University \\[5pt]
 Changchun 130012, P. R. China\\[9pt]
$^2$Center for Combinatorics, LPMC-TJKLC, Nankai University\\[5pt]
 Tianjin 300071, P. R. China\\[9pt]
$^{1,2}$Center for Applied Mathematics, Tianjin University\\[5pt]
 Tianjin 300072, P. R. China\\[9pt]

Email: $^{1}${\tt
qjdou@jlu.edu.cn}, $^{2}${\tt renxy@nankai.edu.cn}
\end{center}

\noindent\textbf{Abstract.} In his study of Ramanujan-Sato type series for $1/\pi$, Sun introduced a sequence of polynomials $S_n(q)$ as given by
$$S_n(q)=\sum\limits_{k=0}^n{n\choose k}{2k\choose k}{2(n-k)\choose n-k}q^k,$$ and he conjectured that the polynomials $S_n(q)$ are $q$-log-convex. By imitating a result of Liu and Wang on generating new $q$-log-convex sequences of polynomials from old ones, we obtain a sufficient condition for determining the $q$-log-convexity of self-reciprocal polynomials. Based on this criterion, we then give an affirmative answer to Sun's conjecture.

\noindent \emph{AMS Classification 2010:} 05A20


\noindent \emph{Corresponding author:} Anne X. Y. Ren, renxy@nankai.edu.cn

\section{Introduction}
The main objective of this paper is to prove a conjecture of Sun \cite{sun2012} on the
$q$-log-convexity of the polynomials $S_n(q)$, which are given by
\begin{align}
S_n(q)&=\sum\limits_{k=0}^n{n\choose k}{2k\choose k}{2(n-k)\choose n-k}q^k\label{conj2}.
\end{align}
These polynomials $S_n(q)$ were introduced by Sun \cite{sun2012} in his study of the Ramanujan-Sato type series for $1/\pi$.

Let us first review some definitions. Recall that a nonnegative sequence $\{a_n\}_{n\geq 0}$ is said to be log-concave if, for any $n\geq 1$,
$$a_n^2\geq a_{n-1}a_{n+1};$$
and is said to be log-convex if, for any $n\geq 1$,
$$a_{n-1}a_{n+1} \geq a_n^2.$$
Many sequences arising in combinatorics, algebra and geometry, turn out to be log-concave or log-convex, see Brenti \cite{bren1994} or Stanley \cite{stan1989}.

For a sequence of polynomials with real coefficients, Stanley introduced the notion of $q$-log-concavity. Throughout we are concerned only with polynomials with real coefficients.
A polynomial sequence $\{f_n(q)\}_{n\geq 0}$ is said to be $q$-log-concave if, for any $n\geq 1$, the difference
$$f_n^2(q)-f_{n+1}(q)f_{n-1}(q)$$
has nonnegative coefficients. The $q$-log-concavity of polynomial sequences has been extensively studied, see Bulter \cite{butl1990}, Krattenthaler \cite{krat1989}, Leroux \cite{lero1990} and Sagan \cite{saga1992}. Similarly, a polynomial sequence $\{f_n(q)\}_{n\geq 0}$ is said to be $q$-log-convex if, for any $n\geq 1$, the difference
$$f_{n+1}(q)f_{n-1}(q)-f_n^2(q)$$
has nonnegative coefficients. Liu and Wang \cite{wang-liu} showed that many classical combinatorial polynomials are $q$-log-convex, see also \cite{chen-yang, chen-wang-yang2010, chen-wang-yang2011}. It should be noted that Butler and Flanigan \cite{butfla2007} introduced a different kind  of $q$-log-convexity.

Sun posed six conjectures on the expansions of $1/\pi$ in terms
of $S_n(q)$, one of which reads
\[\sum_{n=0}^\infty \frac{140n+19}{4624^n}{2n\choose n}S_n(64)=\frac{289}{3\pi}.\]
He also conjectured that the polynomials $S_n(q)$ are $q$-log-convex.
It is easy to see that the coefficients of $S_n(q)$ are symmetric.
Such polynomials are also said to be self-reciprocal. More precisely,
a polynomial
$$f(q)=a_0+a_1q+\cdots+a_nq^n$$
is called a self-reciprocal polynomial of degree $n$ if
$f(q)=q^nf(1/q)$.

In this paper, we shall give a proof of the $q$-log-convexity conjecture of $S_n(q)$. Our proof is closely related to a result of Liu and
Wang, which provides a mechanism of generating new $q$-log-convex sequences of polynomials from certain log-convex sequences of positive numbers and $q$-log-convex sequences of polynomials.
The critical point of their result is to determine the sign of some statistic arising from the difference $f_{n+1}(q)f_{n-1}(q)-f_n^2(q)$ for a given $q$-log-convex sequence $\{f_n(q)\}_{n\geq 0}$.
Assume that $f_n(q)$ has the following form:
\begin{align}
f_n(q)=\sum_{k=0}^n a(n,k)q^k.\label{eq-q-log-convex}
\end{align}
Write the difference $f_{n+1}(q)f_{n-1}(q)-f_n^2(q)$ as
\begin{align*}
\sum_{t=0}^{2n}\left[\sum_{k=0}^{\lfloor t/2\rfloor}\mathcal{\widetilde{L}}_t(a(n,k))
\right]q^t,
\end{align*}
where
\begin{align*}
\mathcal{\widetilde{L}}_t(a(n,k))=\left\{
\begin{array}{ll}
a(n+1,k)a(n-1,t-k)+a(n-1,k)a(n+1,t-k)\\[5pt]
-2a(n,k)a(n,t-k), \qquad \qquad \qquad \mbox{ if }0\leq k< \frac{t}{2},\\[5pt]
a(n+1,k)a(n-1,k)-a^2(n,k), \quad \mbox{ if } t \mbox{ is even and }  k=\frac{t}{2}.
\end{array}
\right.
\end{align*}
Liu and Wang's criterion to determine the $q$-log-convexity of a polynomial sequence is as follows.
\begin{theo}[{\cite[Theorem 4.8]{wang-liu}}]\label{liu-wang thm}
Let $\{u_k\}_{k\geq 0}$ be a log-convex sequence and let $\{f_n(q)\}_{n\geq 0}$ be a $q$-log-convex sequence as defined in \eqref{eq-q-log-convex}. Given $n\geq 1$ and $0\leq t\leq 2n$, if there exists an index $k'$ associated with $n,t$ such that
\begin{align*}
\mathcal{\widetilde{L}}_t(a(n,k))\left\{
\begin{array}{ll}
\geq 0, & \mbox{ if }0\leq k\leq k',\\[5pt]
\leq 0, & \mbox{ if } k'< k\leq \frac{t}{2},
\end{array}
\right.
\end{align*}
then, the polynomial sequence $\{g_n(q)\}_{n\geq 0}$ defined by
\begin{align}\label{eq-gn}
g_n(q)=\sum_{k=0}^n a(n,k)u_kq^k
\end{align}
is $q$-log-convex.
\end{theo}

We attempted to use the above result to prove the $q$-log-convexity of $\{S_n(q)\}_{n\geq 0}$ by taking
\begin{align}\label{seq-ua}
u_k=\binom{2k}{k}, \quad a(n,k)=\binom{n}{k}\binom{2n-2k}{n-k}.
\end{align}
Experimental evidence suggests that $\mathcal{\widetilde{L}}_t(a(n,k))$ meets Liu and Wang's criterion. The determination of the sign of $\mathcal{\widetilde{L}}_t(a(n,k))$ relies on the relative position of two polynomials in $t$ on the interval $[0,2n]$. While it is easier to determine their relative position on the interval $[0, n]$ than on $[0, 2n]$. This forces us to consider the symmetry of the coefficients of the self-reciprocal polynomials to circumvent the above difficulty. As a result, we obtain a criterion for the $q$-log-convexity of self-reciprocal polynomials in the spirit of Theorem \ref{liu-wang thm}, which shall be given in Section 2. By using this criterion, we then confirm the $q$-log-convexity conjecture of Sun in Section 3.

\section{ A criterion for $q$-log-convexity}

The aim of this section is to present a criterion for proving a sequence of self-reciprocal polynomials to be $q$-log-convex.

Noting that for $\mathcal{\widetilde{L}}_t(a(n,t/2))$ in Theorem \ref{liu-wang thm}, only its sign should be considered, we make the following modification to $\mathcal{\widetilde{L}}_t(a(n,k))$ for convenience:
\begin{align}\label{eq-recur-1}
\mathcal{L}_t(a(n,k))=&a(n+1,k)a(n-1,t-k)+a(n-1,k)a(n+1,t-k)\nonumber\\[5pt]
                     &\qquad \qquad \qquad -2a(n,k)a(n,t-k), \qquad \mbox{ if } 0\leq k\leq \frac{t}{2}.
\end{align}
Then we give the following criterion which is applicable to $\{S_n(q)\}_{n\geq 0}$.

\begin{theo}\label{thm-criterion}
Given a log-convex sequence $\{u_k\}_{k\geq 0}$ and a $q$-log-convex sequence $\{f_n(q)\}_{n\geq 0}$ as defined in \eqref{eq-q-log-convex}, let $\{g_n(q)\}_{n\geq 0}$ be the polynomial sequence defined by \eqref{eq-gn}.
Assume that the following two conditions are satisfied:
\begin{itemize}
\item[(C1)]  for each $n\geq 0$, the polynomial $g_n(q)$ is a self-reciprocal polynomial of degree $n$; and

\item[(C2)] for given $n\geq 1$ and $0\leq t\leq n$, there exists an index $k'$ associated with $n,t$ such that \begin{align*}
\mathcal{L}_t(a(n,k))\left\{
\begin{array}{ll}
\geq 0, & \mbox{ if } 0\leq k\leq k',\\[5pt]
\leq 0, & \mbox{ if } k'< k\leq \frac{t}{2}.
\end{array}
\right.
\end{align*}
\end{itemize}
Then, the polynomial sequence $\{g_n(q)\}_{n\geq 0}$
is $q$-log-convex.
\end{theo}

\noindent \textit{Proof of Theorem \ref{thm-criterion}.}
Since each $g_n(q)$ is a self-reciprocal polynomial of degree $n$, we have
\begin{align*}
g_{n-1}(q)&=q^{n-1} g_{n-1}(q^{-1}),\\[5pt]
g_n(q)&=q^n g_n(q^{-1}),\\[5pt]
g_{n+1}(q)&=q^{n+1} g_{n+1}(q^{-1}).
\end{align*}
Therefore,
\begin{align*}
g_n^2(q)&=q^{2n}g_n^2(q^{-1}),\\[5pt]
g_{n-1}(q)g_{n+1}(q)&=q^{2n}g_{n-1}(q^{-1})g_{n+1}(q^{-1}),
\end{align*}
i.e., both $g_{n-1}(q)g_{n+1}(q)$ and $g_n^2(q)$ are self-reciprocal polynomials of degree $2n$.

Writing the difference $g_{n-1}(q)g_{n+1}(q)-g_n^2(q)$ as
$$\sum_{t=0}^{2n} \mathcal{B}(n,t)
q^t,$$
we obtain that, for $0\leq t\leq 2n$,
$$\mathcal{B}(n,t)=\mathcal{B}(n,2n-t)$$
due to reciprocity. Accordingly, to prove the $q$-log-convexity of $\{g_n(q)\}_{n\geq 0}$, it suffices to show that $\mathcal{B}(n,t)$ is nonnegative for any $0\leq t\leq n$.

It is ready to see that
$$\mathcal{B}(n,t)=\left\{
\begin{array}{ll}
\sum_{k=0}^{s}\mathcal{L}_t(a(n,k))u_ku_{t-k}, \mbox{ if $t=2s+1$,} \\[8pt]
\sum_{k=0}^{s-1}\mathcal{L}_t(a(n,k))u_ku_{t-k}+\frac{\mathcal{L}_t\left(a(n,s)\right)}{2}u_s^2, \mbox{ if $t=2s$.}
\end{array}
\right.
$$

To prove the nonnegativity of $\mathcal{B}(n,t)$, we further need to use the log-convexity of $\{u_k\}_{k\geq 0}$ and the $q$-log-convexity of $\{f_n(q)\}_{n\geq 0}$.

On one hand, by the log-convexity of $\{u_k\}_{k\geq 0}$, we have
\begin{align}\label{eq-log-convex}
u_0u_t\geq u_1u_{t-1}\geq \cdots.
\end{align}

On the other hand, if we write
$$f_{n-1}(q)f_{n+1}(q)-f_n^2(q)=\sum_{t=0}^{2n} \mathcal{A}(n,t)
q^t,$$
then
$$\mathcal{A}(n,t)=\left\{
\begin{array}{ll}
\sum_{k=0}^{s}\mathcal{L}_t(a(n,k)), \mbox{ if $t=2s+1$,} \\[8pt]
\sum_{k=0}^{s-1}\mathcal{L}_t(a(n,k))+\frac{\mathcal{L}_t\left(a(n,s)\right)}{2}, \mbox{ if $t=2s$.}
\end{array}
\right.
$$
Since $\{f_n(q)\}_{n\geq 0}$ is $q$-log-convex, for any $0\leq t\leq 2n$, it holds that
$\mathcal{A}(n,t)\geq 0.$

Now we proceed to prove that $\mathcal{B}(n,t)\geq 0$ for $0\leq t\leq n$. We first consider the case when $t$ is odd, namely, $t=2s+1$ for some $s\in \mathbb{N}$.
By \eqref{eq-log-convex} and the condition (C2), we obtain
\begin{align*}
\mathcal{B}(n,t)&=\sum_{k=0}^{s}\mathcal{L}_t(a(n,k))u_ku_{t-k}\\[5pt]
&\geq \sum_{k=0}^{s}\mathcal{L}_t(a(n,k))u_{k'}u_{t-k'}\\[5pt]
&=\mathcal{A}(n,t)u_{k'}u_{t-k'}\geq 0.
\end{align*}
By the same arguments, if $t=2s$ for some $s\in \mathbb{N}$, then
\begin{align*}
\mathcal{B}(n,t)&=\sum_{k=0}^{s-1}\mathcal{L}_t(a(n,k))u_ku_{t-k}
+\frac{\mathcal{L}_t\left(a(n,s)\right)}{2}u_s^2\\[5pt]
&\geq \sum_{k=0}^{s-1}\mathcal{L}_t(a(n,k))u_{k'}u_{t-k'}+\frac{\mathcal{L}_t\left(a(n,s)\right)}{2}u_{k'}u_{t-k'}\\[5pt]
&=\mathcal{A}(n,t)u_{k'}u_{t-k'} \geq 0.
\end{align*}
This completes the proof. \qed

\section{The $q$-log-convexity of $S_n(q)$}

In this section we wish to use Theorem \ref{thm-criterion} to prove Sun's $q$-log-convexity conjecture of $S_n(q)$. The main result of this section is the following theorem.
\begin{theo} \label{thm-conj2}
The polynomials $S_n(q)$ given by \eqref{conj2} form a $q$-log-convex sequence.
\end{theo}

To this end, take
\begin{align}\label{seq-ua}
u_k=\binom{2k}{k}, \quad a(n,k)=\binom{n}{k}\binom{2n-2k}{n-k},
\end{align}
and hence, by \eqref{eq-q-log-convex} and \eqref{eq-gn}, we have
\begin{align}\label{seq-fg}
f_n(q)=\sum\limits_{k=0}^n\binom{n}{k}\binom{2n-2k}{n-k}q^k, \quad g_n(q)=S_n(q).
\end{align}
It is routine to verify that $\{u_k\}_{k\geq 0}$ is a log-convex sequence. It is also clear that $S_n(q)$ is a self-reciprocal polynomial of degree $n$. There remains to show that  the sequence $\{f_n(q)\}_{n\geq 0}$ is $q$-log-convex, and
the triangular array $\{a(n,k)\}_{0\leq k\leq n}$ satisfies the condition (C2) of Theorem \ref{thm-criterion}.
For the former, we have the following result, and for the latter, see Theorem \ref{c1 thm3}.

\begin{theo}\label{thm-q-lqx}
For $n\geq 0$, let $f_n(q)$ be polynomials given by \eqref{seq-fg}.
Then the sequence $\{f_n(q)\}_{n\geq 0}$ is $q$-log-convex.
\end{theo}

\proof
It suffices to show that
the polynomials
$$q^nf_n(q^{-1})=\sum\limits_{k=0}^n\binom{n}{k}\binom{2k}{k}q^k$$
 form a $q$-log-convex sequence. In view of the $q$-log-convexity of $\left\{(1+q)^n\right\}_{n\geq 0}$ and the log-convexity of  $\{\binom{2k}{k}\}_{k\geq 0}$, it is natural to consider whether the triangular array $\{\binom{n}{k}\}_{0\leq k\leq n}$ satisfies the condition of Theorem \ref{liu-wang thm}.

Note that, for $n\geq 1$, $0\leq t\leq 2n$ and $0\leq k\leq t/2$, we have
\begin{align}\label{eqn31-0}
\mathcal{L}_t\left(\binom{n}{k}\right)&=\binom{n+1}{k}\binom{n-1}{t-k}
+\binom{n+1}{t-k}\binom{n-1}{k}-2\binom{n}{t-k}\binom{n}{k}\nonumber\\[5pt]
&=\frac{1}{n(n+1)(n-k+1)}\binom{n}{k}
\binom{n+1}{t-k}\varphi^{(n,t)}(k),
\end{align}
where
\begin{align*}
\varphi^{(n,t)}(x)=&(n+1)(n-x)(n-x+1)+(n+1)(n-t+x)(n-t+x+1)\\[5pt]
&\qquad -2n(n-x+1)(n-t+x+1).
\end{align*}
Thus, the sign of $\mathcal{L}_t\left(\binom{n}{k}\right)$ depends on that of $\varphi^{(n,t)}(k)$ for $0\leq k\leq  t/2$.

To see the sign changes of $\varphi^{(n,t)}(k)$ as $k$ varies from $0$ to $[t/2]$, we consider the values of $\varphi^{(n,t)}(x)$  as $x$ varies over the interval $[0,t/2]$. Taking the derivative of $\varphi^{(n,t)}(x)$ with respect to $x$, we obtain that
$$(\varphi^{(n,t)}(x))'=(4n+2)(2x-t)\leq 0,\ \mbox{for} \ x\leq t/2.$$
Thus $\varphi^{(n,t)}(x)$ is decreasing on the interval $[0,t/2]$. With $\varphi^{(n,t)}(0)=(n+1)(t^2-t)\geq 0$, for given $n$ and $t$, there exists
$k'$ such that
\begin{align*}
\varphi^{(n,t)}(k)\left\{
\begin{array}{ll}
\geq 0, & \mbox{ if }0\leq k\leq k',\\[5pt]
\leq 0, & \mbox{ if } k'< k\leq \frac{t}{2},
\end{array}
\right.
\end{align*}
and hence
\begin{align*}
\mathcal{L}_t\left(\binom{n}{k}\right)\left\{
\begin{array}{ll}
\geq 0, & \mbox{ if }0\leq k\leq k',\\[5pt]
\leq 0, & \mbox{ if } k'< k\leq \frac{t}{2}.
\end{array}
\right.
\end{align*}
By Theorem \ref{liu-wang thm}, we obtain the desired $q$-log-convexity of $\{q^nf_n(q^{-1})\}_{n\geq 0}$. \qed

The remaining part of this section is to prove the following result.

\begin{theo}\label{c1 thm3} Let $\{a(n,k)\}_{0\leq k\leq n}$ be  the triangular array defined by \eqref{seq-ua}. Then, for any $n\geq 1$ and $0\leq t\leq n$, there exists an index $k'$ with respect to $n,t$ such that
\begin{align*}
\mathcal{L}_t(a(n,k))\left\{
\begin{array}{ll}
\geq 0, & \mbox{ if }0\leq k\leq k',\\[5pt]
\leq 0, & \mbox{ if } k'< k\leq \frac{t}{2}.
\end{array}
\right.
\end{align*}
\end{theo}

Before proving Theorem \ref{c1 thm3}, let us make some observations.
For $n\geq 1$, $0\leq t\leq n$ and $0\leq k\leq t/2$, we have
\begin{align*}
\mathcal{L}_t(a(n,k))=&\binom{n+1}{k}\binom{2n-2k+2}{n-k+1}
\binom{n-1}{t-k}\binom{2n-2t+2k-2}{n-t+k-1}\nonumber\\[5pt]
&+\binom{n-1}{k}\binom{2n-2k-2}{n-k-1}\binom{n+1}{t-k}\binom{2n-2t+2k+2}{n-t+k+1}\nonumber\\[5pt]
&-2\binom{n}{k}\binom{2n-2k}{n-k}\binom{n}{t-k}\binom{2n-2t+2k}{n-t+k}.
\end{align*}
By factorization, we obtain
\begin{align}\label{eqn3-2}
\mathcal{L}_t(a(n,k))=&\frac{1}{(n-k+1)^2(n-t+k+1)^2(2n-2k-1)(2n-2t+2k-1)}\nonumber\\[5pt]
                       &\times\frac{1}{n}\binom{n}{k}\binom{2n-2k}{n-k}\binom{n}{t-k}\binom{2n-2t+2k}{n-t+k}\psi^{(n,t)}(k),
\end{align}
where
\begin{align}\label{eqn3-main}
\psi^{(n,t)}(x)=&(n+1)(n-x)^2(n-x+1)^2(2n-2t+2x+1)(2n-2t+2x-1)\nonumber\\[5pt]
       &+(n+1)(n-t+x)^2(n-t+x+1)^2(2n-2x-1)(2n-2x+1)\nonumber\\[5pt]
       &-2n(n-x+1)^2(n-t+x+1)^2(2n-2x-1)(2n-2t+2x-1).
       \end{align}

Clearly, the sign of $\mathcal{L}_t(a(n,k))$ coincides with that of $\psi^{(n,t)}(k)$ unless $t=n$ and $k=0$. Based on this observation, we  divide the proof of Theorem \ref{c1 thm3} into the following three steps:
\begin{itemize}
\item[(S1)] For $n\geq 1$ and $0\leq t\leq n$, prove that $\mathcal{L}_t(a(n,0))\geq 0$ , see Proposition \ref{prop-1};

\item[(S2)] For $n\geq 2$ and $0\leq t\leq n-1$, prove that there exists $k'$ such that
\begin{align*}
\psi^{(n,t)}(k)\left\{
\begin{array}{ll}
\geq 0, & \mbox{ if } 1\leq k\leq k',\\[5pt]
\leq 0, & \mbox{ if } k'< k\leq \frac{t}{2},
\end{array}
\right.
\end{align*}
see Proposition \ref{prop-2};

\item[(S3)] For $n\geq 2$ and $t=n$, prove that there exists $k'$ such that
\begin{align*}
\psi^{(n,n)}(k)\left\{
\begin{array}{ll}
\geq 0, & \mbox{ if } 1\leq k\leq k',\\[5pt]
\leq 0, & \mbox{ if } k'< k\leq \frac{n}{2},
\end{array}
\right.
\end{align*}
see Proposition \ref{prop-3}.
\end{itemize}

Let us first prove the nonnegativity of $\mathcal{L}_t(a(n,0))$.

\begin{prop}\label{prop-1} For any $n\geq 1$ and $0\leq t\leq n$, we have $\mathcal{L}_t(a(n,0))\geq 0$.
\end{prop}
\proof
For $1\leq n\leq 4$, the nonnegativity of $\mathcal{L}_t(a(n,0))$ can be proved directly as follows:
\allowdisplaybreaks
\begin{align*}
\mathcal{L}_0(a(1,0)) &=4,\, \mathcal{L}_1(a(1,0)) =0,\\[5pt]
\mathcal{L}_0(a(2,0)) &=8,\, \mathcal{L}_1(a(2,0)) =8,\, \mathcal{L}_2(a(2,0)) =0,\\[5pt]
\mathcal{L}_0(a(3,0)) &=40,\, \mathcal{L}_1(a(3,0)) =40,\, \mathcal{L}_2(a(3,0)) =46,\, \mathcal{L}_3(a(3,0)) =8,\\[5pt]
\mathcal{L}_0(a(4,0)) &=280,\, \mathcal{L}_1(a(4,0)) =336,\\[5pt]
\mathcal{L}_2(a(4,0)) &=472,\,\mathcal{L}_3(a(4,0)) =332,\, \mathcal{L}_4(a(4,0)) =60.
\end{align*}

For the remainder of the proof, we assume that $n\geq 5$.
It is routine to compute that the sign of $\mathcal{L}_t(a(n,0))$ coincides with that of
\begin{align*}
\frac{\binom{2n}{n}\binom{n}{t}\binom{2n-2t}{n-t}\theta(t)}{n(n+1)(2n-1)(n-t+1)^2(2n-2t-1)},
\end{align*}
where
\begin{align}\label{func-g}
\theta(x)=&(4n^2-1)x^4-2(2n-1)(2n^2+2n+1)x^3+(4n^4+8n^3+8n^2-1)x^2\nonumber\\[5pt]
&-2n(n+1)(2n^2+4n-1)x+2n(2n-1)(n+1)^2.
\end{align}

To prove that $\mathcal{L}_t(a(n,0))\geq 0$,
there are two cases to consider:
\begin{itemize}
\item[(i)] $t=n$. In this case it suffices to show that $\theta(n)\leq 0$. For $n\geq 5$, one can readily check that
\begin{align*}
\theta(n)&=-n(n-1)(n-2)(n+1)< 0.
\end{align*}

\item[(ii)] $0\leq t<n$. In this case it suffices to show that $\theta(t)\geq 0$. To this end, we consider the monotonicity of $\theta(x)$, regarded as a function of $x$, over the interval $[0,n-1]$.
By \eqref{func-g}, taking the derivative of $\theta(x)$ with respect to $x$, we have
\begin{align*}
\theta'(x)&=2(n-x)\theta_1(x),
\end{align*}
where\begin{align*}
\theta_1(x)=&2(1-4n^2)x^2+(2n-1)(2n^2+4n+3)x-(2n^3+6n^2+3n-1).
\end{align*}
We further need the derivative of $\theta_1(x)$:
\begin{equation*}
\theta_1'(x)=(2n-1)\theta_2(x),
\end{equation*}
where
$$\theta_2(x)=-4(2n+1)x+(2n^2+4n+3).$$
Note that, for $n\geq 5$,
\begin{align*}
\theta_2(0)&=2n^2+4n+3>0,\quad \theta_2(n-1)=-6n^2+8n+7<0.
\end{align*}
Therefore, $\theta_2(x)$ decreases from a positive value to a negative value as $x$ increases from $0$ to $n-1$. Hence, $\theta_1(x)$ first increases and then decreases as $x$ increases from $0$ to $n-1$.
Since, for $n\geq 5$,
\begin{align*}
\theta_1(0)&=1-2n^3-6n^2-3n<0,\\[5pt]
\theta_1(1)&=n(2(n-2)^2-9)>0,\\[5pt]
\theta_1(n-1)&=-4n^4+16n^3-16n^2-12n+6<0,
\end{align*}
there exist $0< x_1<x_2< n-1$ such that
\begin{align*}
\theta_1(x)\left\{
\begin{array}{ll}
< 0, & \mbox{ if } x\in [0,x_1),\\[5pt]
\geq 0, & \mbox{ if } x\in [x_1,x_2],\\[5pt]
< 0, & \mbox{ if } x\in (x_2,n-1].
\end{array}
\right.
\end{align*}
Thus, $\theta(x)$ is decreasing on the interval $[0,x_1)$, increasing on $[x_1,x_2]$, and decreasing on $(x_2,n-1]$.

Note that, for $n\geq 5$, we have
\begin{align*}
\theta(0)&=2n(2n-1)(n+1)^2>0,\\[5pt]
\theta(1)&=2n^2(2n-1)(n-1)>0,\\[5pt]
\theta(2)&=2(n-2)(6n^3-13n^2+1)>0,\\[5pt]
\theta(n-1)&=-4+8n+3n^4-10n^3+11n^2>0.
\end{align*}
It is easy to check that
$$\theta(0)>\theta(1)<\theta(2)>\theta(n-1).$$
By virtue of the monotonicity of $\theta(x)$ on the interval $[0,n-1]$, we must have $x_1\leq 2$. If $x_2> 2$, then $\theta(x)$ is increasing
on $[2,x_2]$, and decreasing on $(x_2,n-1]$. If $x_2\leq 2$, then $\theta(x)$ decreases on $(2,n-1]$. In both cases, we obtain that $\theta(x)>0$ for $x\in[2,n-1]$. In view of $\theta(0)>0$ and $\theta(1)>0$, it is clear that $\theta(t)>0$ for any integer $0\leq t\leq n-1$.
\end{itemize}
Combining (i) and (ii), we obtain the desired result.
\qed

Now we proceed to determine the sign of $\psi^{(n,t)}(k)$ for $n\geq 2$ and $0\leq t\leq n-1$.

\begin{prop}\label{prop-2}
Given $n\geq 2$ and $0\leq t\leq n-1$,
there exists $k'$ with respect to $n,t$ such that
\begin{align*}
\psi^{(n,t)}(k)\left\{
\begin{array}{ll}
\geq 0, & \mbox{ if } 1\leq k\leq k',\\[5pt]
\leq 0, & \mbox{ if } k'< k\leq \frac{t}{2}.
\end{array}
\right.
\end{align*}
\end{prop}

\proof By \eqref{eqn3-2} and Proposition \ref{prop-1}, we know that
$\psi^{(n,t)}(0)\geq 0$. Therefore, it suffices to prove that there exists $0\leq t_0\leq t/2$ such that $\psi^{(n,t)}(x)$, regarded as a function of $x$, is increasing on the interval $[0,t_0)$ and decreasing on the interval $[t_0,t/2]$. To this end, we need to determine the sign changes of the derivative of $\psi^{(n,t)}(x)$ with respect to $x$ on the interval $[0,t/2]$.

Taking the derivative of $\psi^{(n,t)}(x)$, we obtain that
\begin{align*}
(\psi^{(n,t)}(x))'&=2(2x-t)\psi^{(n,t)}_1(x),
\end{align*}
where
\allowdisplaybreaks
\begin{align*}
\psi^{(n,t)}_1(x)=&12(2n+1)x^4-24t(2n+1)x^3\\[5pt]
&-2(16n^3-8(2t-1)n^2-2(7t^2+3t+1)n-(8t^2-4t+3))x^2\\[5pt]
&+2t(16n^3-8(2t-1)n^2-2(t^2+3t+1)n-(2t^2-4t+3))x\\[5pt]
&+\left(8n^5-4(4t-1)n^4+4(t^2-t-3)n^3+4(-t^2+5t+t^3-2)n^2\right.\\[5pt]
&\left.+(4t^3-10t^2-1+11t)n-(2t^2-3t+1)\right).
\end{align*}
We further need to consider the derivative of $\psi^{(n,t)}_1(x)$:
\begin{align}\label{eq-diffh1}
(\psi^{(n,t)}_1(x))'&=2(2x-t)\psi^{(n,t)}_2(x),
\end{align}
where
\allowdisplaybreaks
\begin{align*}
\psi^{(n,t)}_2(x)=&12(2n+1)x^2-12t(2n+1)x-16n^3+8(2t-1)n^2\\[5pt]
& +2(t^2+3t+1)n+(2t^2-4t+3).
\end{align*}
Note that the axis of symmetry of the quadratic function $\psi^{(n,t)}_2(x)$ is $x=t/2$. Hence, $\psi^{(n,t)}_2(x)$ decreases as $x$ increases from $0$ to $t/2$. It is routine to verify that, for $n\geq 1$ and $0\leq t<n$,
\begin{align*}
\psi^{(n,t)}_2\left(\frac{t}{2}\right)=-4n(2n-t)^2-(4n-t-1)(2n-t)-3(t-1)<0.
\end{align*}
Let $x_0$ be the zero of $\psi^{(n,t)}_2(x)$ to the left of the axis of symmetry. Then we have
\begin{align*}
\psi^{(n,t)}_2(x)\left\{
\begin{array}{ll}
> 0,& \mbox{ if } 0\leq x< x_0,\\[5pt]
< 0,& \mbox{ if } x_0< x< t/2.
\end{array}
\right.
\end{align*}
By \eqref{eq-diffh1}, we have
\begin{align*}
(\psi^{(n,t)}_1(x))'\left\{
\begin{array}{ll}
< 0,& \mbox{ if } 0\leq x< x_0,\\[5pt]
> 0,& \mbox{ if } x_0<x< t/2.
\end{array}
\right.
\end{align*}
If $x_0\leq 0$, this means that $\psi^{(n,t)}_1(x)$ is increasing on $[0,t/2]$. If $x_0> 0$, this means that $\psi^{(n,t)}_1(x)$ is decreasing on $[0,x_0]$ and increasing on $[x_0,t/2]$.

We proceed to determine the sign changes of $(\psi^{(n,t)}(x))'$ based on the  above monotonicity of $\psi^{(n,t)}_1(x)$. For our purpose, the values of $\psi^{(n,t)}_1(x)$ at the two endpoints of the interval $[0,t/2]$ are to be examined.

We claim that, for any integers $n\geq 2$ and $0\leq t<n$,
it holds $\psi^{(n,t)}_1(t/2)>0$. Using Maple, we find that
\begin{align*}
\psi^{(n,t)}_1\left(\frac{t}{2}\right) =&8n^5-16n^4t+12n^3t^2-4n^2t^3+\frac{1}{2}nt^4+4n^4-4n^3t+nt^3-\frac{1}{4}t^4\\[5pt]
&-12n^3+20n^2t-11nt^2+2t^3-8n^2+11nt-\frac{7}{2}t^2-n+3t-1\\[5pt]
=&\left(\frac{1}{2}n-\frac{1}{4}\right)(2n-t)^4+(n-2)(2n-t)^3+\left(n-\frac{7}{2}\right)(2n-t)^2\\[5pt]
&+3(n-1)(2n-t)+5n-1,
\end{align*}
which is greater than $0$ whenever $n\geq 4$ and $0\leq t<n$. It remains to check the validity of $\psi^{(n,t)}_1(t/2)>0$ for $n=2,3$. In fact,  for $n=2$, we have $0\leq t<2$ and hence
\begin{align*}
\psi^{(2,t)}_1\left(\frac{t}{2}\right)=
\frac{3}{4}\left((4-t)^2-1\right)^2+3(4-t)+\frac{33}{4}
>0.
\end{align*}
For $n=3$, we have $0\leq t<3$ and hence
\begin{align*}
\psi^{(3,t)}_1\left(\frac{t}{2}\right)=
\frac{5}{4}(6-t)^4+\left(\frac{11}{2}-t\right)(6-t)^2+6(6-t)+14>0.
\end{align*}

As we see, the value of $\psi^{(n,t)}_1(t/2)$ must be positive. By further taking into account the value of $x_0$ and the sign of $\psi^{(n,t)}_1(0)$, there are three cases to determine the monotonicity of $\psi^{(n,t)}(x)$:

\begin{itemize}
\item[(i)] $x_0\leq 0$ and $\psi^{(n,t)}_1(0)\geq 0$. In this case, $\psi^{(n,t)}_1(x)$ increases from a nonnegative value to a positive value as $x$ increases from $0$ to $t/2$. Thus, $(\psi^{(n,t)}(x))'$ takes only nonpositive values on $[0,t/2]$. That is to say, $\psi^{(n,t)}(x)$ is decreasing on the interval $[0,t/2]$.

\item[(ii)] $x_0\leq 0$ and $\psi^{(n,t)}_1(0)<0$. In this case, $\psi^{(n,t)}_1(x)$ increases from a negative value to a positive value as $x$ increases from $0$ to $t/2$. Therefore, there exists $0<t_0<t/2$ such that
    \begin{align*}
      \psi^{(n,t)}_1(x)\left\{
      \begin{array}{ll}
         \leq 0,& \mbox{ if } 0\leq x\leq t_0,\\[5pt]
         \geq 0,& \mbox{ if } t_0<x\leq t/2.
      \end{array}
      \right.
    \end{align*}
Hence,   we have
  \begin{align*}
      (\psi^{(n,t)}(x))'\left\{
      \begin{array}{ll}
         \geq 0,& \mbox{ if } 0\leq x\leq t_0,\\[5pt]
         \leq 0,& \mbox{ if } t_0<x\leq t/2.
      \end{array}
      \right.
    \end{align*}
That is to say, $\psi^{(n,t)}(x)$ is increasing on $[0,t_0]$ and decreasing on $[t_0,t/2]$.

\item[(iii)] $0< x_0<t/2$. In this case, we must have $\psi^{(n,t)}_1(0)< 0$. Once this assertion is proved, we obtain the desired monotonicity of $\psi^{(n,t)}(x)$ on $[0,t/2]$, by using similar arguments as in case (ii).
    Note that the condition $0< x_0<t/2$ implies that $\psi^{(n,t)}_2(0)>0$.

    Now we are to deduce $\psi^{(n,t)}_1(0)< 0$ from the positivity of $\psi^{(n,t)}_2(0)$. Using maple, we find that
    \begin{align*}
    \psi^{(n,t)}_1(0)=&(n+1)\left(4nt^3+2(2n^2-4n-1)t^2-(16n^3-12n^2-8n-3)t\right.\\
           &\left.\qquad\qquad\qquad+(8n^4-4n^3-8n^2-1)\right),\\[5pt]
    \psi^{(n,t)}_2(0)=&2(n+1)t^2+2(8n^2+3n-2)t-(2n-1)(8n^2+8n+3).
    \end{align*}

    Recall that $0\leq t\leq n-1$ by the hypothesis. We may regard $\psi^{(n,t)}_1(0)/(n+1)$ as a polynomial in the variable $t$ over the interval $[0,n-1]$, denoted by $\xi(t)$, and similarly, regard $\psi^{(n,t)}_2(0)$ as a polynomial $\eta(t)$. Now we can divide the proof of $\psi^{(n,t)}_1(0)< 0$ into the following three statements:

    \medskip

    \textbf{Claim 1.} If $\psi^{(n,t)}_2(0)>0$, then $n\neq 2,3$.

    \noindent \textit{Proof of Claim 1.} In fact, it is routine to check that $\psi^{(n,t)}_2(0)<0$ if $(n,t)\in \{(2,0),(2,1),(3,0),(3,1),(3,2)\}$, contradicting the positivity of $\psi^{(n,t)}_2(0)$. \qed

    \medskip

    \noindent \textbf{Claim 2.} For any integer $n\geq 4$, the polynomial $\xi(t)$ takes only negative values on the interval $[\frac{3}{4}n,n-1]$.

    \noindent \textit{Proof of Claim 2.} Note that, for $n\geq 4$, it is routine to check that
\begin{align*}
\xi\left(\frac{3}{4}n\right)&=-\frac{1}{64}\left(4n^2(n-4)^2+136\left(n-\frac{9}{17}\right)^2+\frac{440}{17}\right)<0,\\
\xi(n-1)&=-(4n-18)n^2-13n-6<0.
\end{align*}
We further need to consider the first order derivative and the second order derivative of $\xi(t)$ with respect to $t$:
\begin{align*}
\xi'(t)&=12nt^2+(8n^2-16n-4)t+(12n^2-16n^3+8n+3),\\[5pt]
\xi''(t)&=24nt+(8n^2-16n-4).
\end{align*}
Note that, for $n\geq 4$, we have
$$
\xi''(0)=8(n-1)^2-12>0.
$$
Hence, $\xi''(t)>0$ for any $0\leq t\leq n-1$. That is to say, $\xi'(t)$ is strictly increasing on the interval $[0,n-1]$. It is clear that, for $n\geq 4$,
$$\xi'\left(\frac{3}{4}n\right)=-\frac{13}{4}n^3+5n+3<0.$$
Thus, there exists $3n/4\leq t_1\leq n-1$ such that
    \begin{align*}
      \xi'(t)\left\{
      \begin{array}{ll}
         \leq 0,& \mbox{ if } \frac{3}{4}n\leq t\leq t_1,\\[5pt]
         > 0,& \mbox{ if } t_1<t\leq n-1.
      \end{array}
      \right.
    \end{align*}
In view of $\xi(3n/4)<0$ and $\xi(n-1)<0$, we obtain $\xi(t)<0$ for any $t\in[\frac{3}{4}n,n-1]$. This completes the proof of Claim 2.\qed

\medskip

\noindent \textbf{Claim 3.} For any integer $n\geq 2$, the polynomial $\eta(t)$ takes only negative values on the interval $[0,\frac{3}{4}n]$.

\noindent \textit{Proof of Claim 3.} For $n\geq 2$, a straightforward computation shows that
\begin{align*}
\eta(0)&=-16n^3-8n^2+2n+3<0,\\[5pt]
\eta\left(\frac{3}{4}n\right)&=-\frac{23}{8}n^3-\frac{19}{8}n^2-n+3<0.
\end{align*}
Note that the axis of symmetry of the quadratic function $\eta(t)$ is $$t=-\frac{8n^2+3n-2}{2(n+1)},$$
which lies strictly to the left of $y$-axis. Therefore, $\eta(t)<0$ for any $t\in [0,3n/4]$ since both $\eta(0)$ and $\eta(3n/4)$ are negative. This completes the proof of Claim 3.\qed

\medskip
Now we can prove the negativity of $\psi^{(n,t)}_1(0)$. From   $\psi^{(n,t)}_2(0)>0$ it follows $\eta(t)>0$. By Claim 3, we must have $t> 3n/4$. Then by Claim 1 and Claim 2, we get $\xi(t)<0$, and hence $\psi^{(n,t)}_1(0)<0$, as desired.
\end{itemize}
Combining (i), (ii) and (iii), we complete the proof.
\qed

The above proposition is the key step for the proof of  Theorem \ref{c1 thm3}. Finally, we need to establish the following result.

\begin{prop}\label{prop-3}
Given $n\geq 2$, there exists $k'$ with respect to $n$ such that
\begin{align*}
\psi^{(n,n)}(k)\left\{
\begin{array}{ll}
\geq 0, & \mbox{ if } 1\leq k\leq k',\\[5pt]
\leq 0, & \mbox{ if } k'< k\leq \frac{n}{2}.
\end{array}
\right.
\end{align*}
\end{prop}

\proof  By \eqref{eqn3-main}, we obtain that
\begin{align*}
\psi^{(n,n)}(x)=&8(2n+1)x^6-24n(2n+1)x^5+2(26n^3-2n+12n^2+3)x^4\\[5pt]
    &-4n(3+6n^3+2n^2-2n)x^3+2(4n^2+2n-1-4n^3+2n^5)x^2\\[5pt]
    &+2n(n-1)(2n-1)(n+1)x-n(n-1)(n-2)(n+1)^2.
\end{align*}

It is easy to check that $\psi^{(2,2)}(1)=8$. Hence, the proposition holds for $n=2$. For the remainder of the proof, assume that $n\geq 3$. To determine the sign of $\psi^{(n,n)}(k)$, let us consider the derivative of $\psi^{(n,n)}(x)$ with respect to $x$. Using Maple, we get
\begin{align*}
(\psi^{(n,n)}(x))'&=2(2x-n)\psi^{(n,n)}_1(x),
\end{align*}
where
\begin{align*}
\psi^{(n,n)}_1(x)=&12(1+2n)x^4-24n(1+2n)x^3+2(6n^2-2n+3+14n^3)x^2\\[5pt]
&-2n(2n^3+3-2n)x-(n-1)(2n-1)(n+1).
\end{align*}
We also need to consider the derivative of $\psi^{(n,n)}_1(x)$ with respect to $x$:
\begin{align*}
(\psi^{(n,n)}_1(x))'&=2(2x-n)\psi^{(n,n)}_2(x),
\end{align*}
where
\begin{align*}
\psi^{(n,n)}_2(x)=12(1+2n)x^2-12n(1+2n)x+2n^3+3-2n.
\end{align*}
Note that the axis of symmetry of the quadratic function $\psi^{(n,n)}_2(x)$ is $x=n/2$, and, for $n\geq 3$,
\begin{align*}
\psi^{(n,n)}_2(0)&=2n^3-2n+3>0,\\[5pt]
\psi^{(n,n)}_2(n/2)&=-4n^3-3n^2-2n+3<0.
\end{align*}
Thus, $\psi^{(n,n)}_2(x)$ decreases from a positive value to a negative value as $k$ increases from $0$ to ${n}/{2}$. Hence, there exists $0<x_0<n/2$ such that
    \begin{align*}
      (\psi^{(n,n)}_1(x))'\left\{
      \begin{array}{ll}
         \leq 0,& \mbox{ if } 0\leq x\leq x_0,\\[5pt]
         \geq 0,& \mbox{ if } x_0<x\leq n/2.
      \end{array}
      \right.
    \end{align*}
In view of that, for $n\geq 3$,
\begin{align*}
\psi^{(n,n)}_1(0)&=-n^2(n-1)-n(n^2-2)-1<0,\\[5pt]
\psi^{(n,n)}_1(n/2)&=\frac{1}{4}(2n^3(n^2-2)+n^2(3n^2-2)+4(2n-1))>0,
\end{align*}
there exists $0<x_1<n/2$ such that
    \begin{align*}
      \psi^{(n,n)}_1(x)\left\{
      \begin{array}{ll}
         \leq 0,& \mbox{ if } 0\leq x\leq x_1,\\[5pt]
         \geq 0,& \mbox{ if } x_1<x\leq n/2.
      \end{array}
      \right.
    \end{align*}
Therefore,
    \begin{align*}
      (\psi^{(n,n)}(x))'\left\{
      \begin{array}{ll}
         \geq 0,& \mbox{ if } 0\leq x\leq x_1,\\[5pt]
         \leq 0,& \mbox{ if } x_1<x\leq n/2.
      \end{array}
      \right.
    \end{align*}
Moreover, it is easy to verify that, for $n\geq 3$,
\begin{align*}
\psi^{(n,n)}(1)&=(n-1)((3n-16)n^3+(21n^2+8n-12))>0, \\[5pt]
\psi^{(n,n)}(n/2)&=-\frac{1}{8}n(n-1)(n^2-n-4)(n+2)^2<0.
\end{align*}
Thus, there exists $1<x_2<n/2$ such that
    \begin{align*}
      \psi^{(n,n)}(x)\left\{
      \begin{array}{ll}
         \geq 0,& \mbox{ if } 1\leq x\leq x_2,\\[5pt]
         \leq 0,& \mbox{ if } x_2<x\leq n/2.
      \end{array}
      \right.
    \end{align*}
Thus, there exists an index $k'=k'(n,n)$ such that $\psi^{(n,n)}(k)\geq 0$ for $1\leq k\leq k'$ and $\psi^{(n,n)}(k)\leq 0$ for $k'<k\leq n/2$, as desired. This completes the proof.\qed

We now come to the proof of Theorem \ref{c1 thm3}.

\noindent \textit{Proof of Theorem \ref{c1 thm3}.} By Proposition \ref{prop-1}, for any $n\geq 1$ and $0\leq t\leq n$, we have $\mathcal{L}_t(a(n,0))\geq 0$. Given $n\geq 1$, it suffices to show that, for $0\leq t\leq n$, there exists $k'$ such that $\mathcal{L}_t(a(n,k))\geq 0$ for $1\leq k\leq k'$ and $\mathcal{L}_t(a(n,k))\leq 0$ for $k'< k\leq t/2$. By \eqref{eqn3-2}, for $k\geq 1$, the sign of $\mathcal{L}_t(a(n,k))$ coincides with that of $\psi^{(n,t)}(k)$. Combining Propositions \ref{prop-2} and \ref{prop-3}, we obtain the desired result. \qed

Finally, we can prove the $q$-log-convexity of $\{S_n(q)\}_{n\geq 0}$.

\noindent \textit{Proof of Theorem \ref{thm-conj2}.} This immediately follows from Theorems \ref{thm-criterion}, \ref{thm-q-lqx} and \ref{c1 thm3}.

\vskip 3mm
\noindent {\bf Acknowledgments.} This work was supported by the 973 Project, the PCSIRT Project
of the Ministry of Education, the National Science Foundation of China, and the Fundamental Research Funds for the Central Universities of China.

\end{document}